\documentclass{amsart}
\usepackage{mathrsfs}
\usepackage{mathrsfs}
\usepackage{amssymb}
\usepackage{amsmath,amssymb}

\newtheorem{theorem}{Theorem}[section]
\newtheorem{lemma}[theorem]{Lemma}
\newtheorem{proposition}[theorem]{Proposition}

\theoremstyle{definition}

\theoremstyle{remark}

\numberwithin{equation}{section}

\begin{document}

\title[Several analytic  inequalities] {Several analytic  inequalities  in some $Q-$spaces}

\author{Pengtao Li}
\address{School of Mathematics, Peking University, Beijing, 100871, China}
\curraddr{Department of Mathematics and Statistics, Memorial
University of Newfoundland, St. John's, NL A1C 5S7, Canada}

\author{Zhichun Zhai}
\address{Department of Mathematics and Statistics, Memorial University of Newfoundland, St. John's, NL A1C 5S7, Canada}
\curraddr{}
\thanks{Project supported in part  by Natural Science and
Engineering Research Council of Canada.}
\thanks{Corresponding author: Zhichun Zhai}
\thanks{E-mail addresses: li\_ptao@163.com (Pengtao Li);
a64zz@mun.ca(Zhichun Zhai)}


\date{}


\begin{abstract} In this paper,  we establish separate necessary and sufficient
  John-Nirenberg (JN) type inequalities  for functions in
 $Q_{\alpha}^{\beta}(\mathbb{R}^{n})$ which  imply Gagliardo-Nirenberg (GN) type inequalities in
 $Q_{\alpha}(\mathbb{R}^{n}).$ Consequently,  we obtain Trudinger-Moser
 type inequalities  and Brezis-Gallouet-Wainger type inequalities in $Q_{\alpha}(\mathbb{R}^{n}).$
  \end{abstract}

\maketitle


 \vspace{0.1in}

 \section{Introduction and Statement of Main Results}
  This paper  studies  several analytic  inequalities in some $Q$ spaces. We first establish
   John-Nirenberg  type inequalities  in $Q_{\alpha}^{\beta}(\mathbb{R}^{n})(n\geq
   2).$  Then we get  Gagliardo-Nirenberg,   Trudinger-Moser
  and Brezis-Gallouet-Wainger type inequalities
in $Q_{\alpha}(\mathbb{R}^{n}).$
 Here $Q_{\alpha}^{\beta}(\mathbb{R}^{n})$ is   the set of all measurable complex-valued functions $f$ on
 $\mathbb{R}^{n}$  satisfying
 \begin{equation}\label{def Q a b}
\|f\|_{Q_{\alpha}^{\beta}(\mathbb{R}^{n})}
 =\sup_{I}\left((l(I))^{2(\alpha+\beta-1)-n}\int_{I} \int_{I}\frac{|f(x)-f(y)|^{2}}{|x-y|^{n+2(\alpha-\beta+1)}}dxdy\right)^{1/2}<\infty
\end{equation}
 for  $\alpha\in (-\infty,\beta)$ and $\beta\in (1/2,1],$ where the supremum is taken over all cubes $I$ with the edge length $l(I)$
 and the edges parallel to the coordinate axes in $\mathbb{R}^{n}.$ Obviously, $Q_{\alpha}^{1}(\mathbb{R}^{n})=Q_{\alpha}(\mathbb{R}^{n})$
 which was introduced by Essen, Janson, Peng and Xiao in  \cite{M Essen S Janson L Peng J Xiao}. It has been found that
 $Q_{\alpha}(\mathbb{R}^{n})$ is a useful and interesting concept, see, for example, Dafni and Xiao \cite{G. Dafni J. Xiao, G Dafni J Xiao 1},
  Xiao \cite{J. Xiao 1}, Cui and Yang \cite{L Cui Q Yang}. As a generalization of $Q_{\alpha}(\mathbb{R}^{n}),$ $Q_{\alpha}^{\beta}(\mathbb{R}^{n})$
  is very useful in harmonic analysis and partial differential equations, see  Yang and Yuan \cite{D. Yang W. Yuan}, Li and Zhai
  \cite{Pengtao Li Zhicun Zhai, Pengtao Li Zhichun Zhai 1} in which $Q_{\alpha}^{\beta}(\mathbb{R}^{n})$ was applied to study the
 well-posednes and regularity of mild solutions to fractional Navier-Stokes equations with fractional Laplacian $(-\triangle)^{\beta}.$

 JN type  inequality is  classical  in modern analysis and widely applied in theory of partial differential equations. In \cite{John Nirenberg}, John and Nirenberg
 proved the JN inequality for $BMO(\mathbb{R}^{n}).$
 In this paper,     we establish JN type inequalities in $
 Q_{\alpha}^{\beta}(\mathbb{R}^{n})$ a special case of which
implies  Gagliardo-Nirenberg (GN) type inequalities meaning  the
continuous embeddings such as $L^{r}(\mathbb{R}^{n})\cap
Q_{\alpha}(\mathbb{R}^{n})\subseteq
 L^{p}(\mathbb{R}^{n})$ for $-\infty<\alpha<1$ and $1\leq r\leq p<\infty.$
  Moreover,
  from GN type inequalities in $Q_{\alpha}(\mathbb{R}^{n}),$ we get Trudinger-Moser  and Brezis-Gallouet-Wainger type
  inequalities. See, for example,
  \cite{Brezis Gallouet, Brezis Wainger, Engler,  Kozono H. Sato T. Wadade H.,  Kozono H. Taniuchi Y.} for more information  about
 Trudinger-Moser and Brezis-Gallouet-Wainger type inequalities. To achieve our main goals, we  need  the characterization of $Q_{\alpha}^{\beta}(\mathbb{R}^{n})$
 in terms of the square mean oscillation over cubes.

 We recall some facts about   mean oscillation over
 cubes. For any cube $I$ and an integrable function $f$ on $I,$ we define
\begin{equation}\label{mean ch Q a b 1}
f(I)=\frac{1}{|I|}\int_{I}f(x)dx
\end{equation}
the mean of $f$ on $I,$ and for $1\leq q<\infty,$
 \begin{eqnarray}\label{mean ch Q a b 2}
  \Phi_{f}^{q}(I)=\frac{1}{|I|}\int_{I}|f(x)-f(I)|^{q}dx
 \end{eqnarray}
 the $q-$mean oscillation of $f$ on $I.$ Recall  the well-known identities
\begin{equation}\label{mean ch Q a b 3}
\frac{1}{|I|}\int_{I}|f(x)-a|^{2}dx=\Phi_{f}^{2}(I)+|f(I)-a|^{2}
\end{equation}
for any complex number $a,$ and
\begin{equation}\label{mean ch Q a b 4}
\frac{1}{|I|^{2}}\int_{I}\int_{I}|f(x)-f(y)|^{2}dxdy=2\Phi^{2}_{f}(I).
\end{equation}
Moreover, if $I\subset J,$ then we have
\begin{equation}\label{mean ch Q a b 5}
 \Phi^{2}_{f}(I)\leq \frac{|J|}{|I|}\Phi^{2}_{f}(J)
\end{equation}
and
\begin{equation}\label{mean ch Q a b 6}
|f(I)-f(J)|^{2}\leq \frac{|J|}{|I|}\Phi^{2}_{f}(J).
\end{equation}

 Let $\mathcal {D}_{0}=\mathcal{D}_{0}(\mathbb{R}^{n})$ be the set of unit cubes whose vertices have integer coordinates, and let, for
  any integer $k\in \mathbb{Z},$ $\mathcal{D}_{k}=\mathcal{D}_{k}(\mathbb{R}^{n})=\{2^{-k}I:I\in \mathcal{D}_{0}\},$  then the cubes in
  $\mathcal{D}=\cup_{-\infty}^{\infty}\mathcal{D}_{k}$ are called dyadic. Furthermore, if $I$ is any cube,  $\mathcal{D}_{k}(I),$ $k\geq 0,$
  denote    the set of the $2^{kn}$  subcubes of edge length $2^{-k}l(I)$ obtained by $k$ successive bipartitions of each edge of $I.$ Moreover,
  put $\mathcal{D}(I)=\cup_{0}^{\infty}\mathcal{D}_{k}(I).$ For any cube $I$ and a measurable function $f$ on $I,$ we define
\begin{eqnarray}\label{mean ch Q a b 7}
\Psi_{f,\alpha,\beta}(I)&=&(l(I))^{4\beta-4}\sum_{k=0}^{\infty}\sum_{J\in
\mathcal{D}_{k}(I)}2^{(2(\alpha-\beta+1)-n)k}\Phi^{2}_{f}(J)\\
&=&(l(I))^{4\beta-4}\sum_{J\in
\mathcal{D}(I)}\left(\frac{l(J)}{l(I)}\right)^{n-2(\alpha-\beta+1)}\Phi^{2}_{f}(J)\nonumber.
\end{eqnarray}

 We can prove the following  proposition by a similar argument   applied by   Essen, Janson, Peng and Xiao for the case $\beta=1$ in
 \cite[Theorem 5.5]{M Essen S Janson L Peng J Xiao}. The details are omitted  here.

 \begin{proposition}\label{mean ch Q a b theorem1}
 Let $-\infty<\alpha<\beta$ and $\beta\in (1/2,1].$ Then $Q_{\alpha}^{\beta}(\mathbb{R}^{n})$ equals  the space of all  measurable functions
  $f$ on $\mathbb{R}^{n}$ such that $\sup_{I}\Psi_{f,\alpha,\beta}(I)$ is finite, where $I$ ranges over  all cubes in $\mathbb{R}^{n}.$ Moreover,
  the square root of this supremum is a norm on $Q_{\alpha}^{\beta}(\mathbb{R}^{n}),$ equivalent to $\|f\|_{Q_{\alpha}^{\beta}(\mathbb{R}^{n})}$ as
 defined above.
 \end{proposition}

  Using this equivalent characterization of $Q_{\alpha}^{\beta}(\mathbb{R}^{n}),$ we can establish the
  following JN type inequalities.
\begin{theorem}\label{mean ch Q a b theorem2}
 Let $-\infty< \alpha<\beta,$  $\beta\in (1/2,1]$ and  $0\leq p<2.$ If there exist positive constants $B,C$ and $c,$ such that, for all
cubes $I\subset \mathbb{R}^{n},$ and any $t>0,$
\begin{equation}\label{mean ch Q a b 113}
 (l(I))^{4\beta-4}\sum_{k=0}^{\infty}2^{(2(\alpha-\beta+1)-n)k}\sum_{J\in \mathcal{D}_{k}(I)}\frac{m_{J}(t)}{|J|}\leq
B\max\left\{1,\left(\frac{C}{t}\right)^{p}\right\}\exp(-ct),
\end{equation}
 then $f$ is a function in $Q_{\alpha}^{\beta}(\mathbb{R}^{n}).$ Here $m_{I}(t)$  is the distribution function of $f-f(I)$ on the cube $I:$
\begin{equation}\label{mean ch Q a b 13}
m_{I}(t)=|\{x\in I: |f(x)-f(I)|>t\}|.
\end{equation}
\end{theorem}

\begin{theorem}\label{mean ch Q a b theorem3}
 Let $-\infty< \alpha<\beta,$  $\beta\in (1/2,1]$  and $f\in Q_{\alpha}^{\beta}(\mathbb{R}^{n}).$ Then there exist positive
constants $B$ and $b,$ such that
\begin{equation}\label{mean ch Q a b 14}
(l(I))^{4\beta-4}\sum_{k=0}^{\infty}2^{(2(\alpha-\beta+1)-n)k}\sum_{J\in
\mathcal{D}_{k}(I)}\frac{m_{J}(t)}{|J|}\leq
B\max\left\{1,\left(\frac{\|f\|_{Q_{\alpha}^{\beta}}}{t}\right)^{2}\right\}\exp\left(\frac{-bt}{\|f\|_{Q_{\alpha}^{\beta}}}\right)
\end{equation}
 holds for  $t\leq \|f\|_{Q_{\alpha}^{\beta}(\mathbb{R}^{n})}$ and any cubes $I\subset \mathbb{R}^{n},$ or
 for $t>\|f\|_{Q_{\alpha}^{\beta}(\mathbb{R}^{n})}$ and   cubes $I\subset \mathbb{R}^{n}$ with $(l(I))^{2\beta-2}\geq1.$ Moreover,
   there  holds
\begin{equation}\label{mean ch Q a b 14''}
(l(I))^{4\beta-4}\sum_{k=0}^{\infty}2^{(2(\alpha-\beta+1)-n)k}\sum_{J\in
\mathcal{D}_{k}(I)}\frac{m_{J}(t)}{|J|}\leq B
\end{equation}
 for  $t>\|f\|_{Q_{\alpha}^{\beta}(\mathbb{R}^{n})}$ and   cubes $I\subset \mathbb{R}^{n}$ with $(l(I))^{2\beta-2}<1.$
\end{theorem}

 For $\beta=1,$    the JN inequality in $Q_{\alpha}(\mathbb{R}^{n})$ was conjectured by Essen-Janson-Peng-Xiao in
  \cite{M Essen S Janson L Peng J Xiao} and finally a modified version as in Theorems \ref{mean ch Q a b theorem2}-\ref{mean ch Q a b theorem3}
   was established  by Yue-Dafni \cite{H Yue G Dafni}.

 According to  Essen, Janson, Peng and  Xiao \cite[Theorem 2.3]{M Essen S Janson L Peng J Xiao} and Li and Zhai \cite[Theorem 3.2]{Pengtao Li Zhicun Zhai}, we know that if $-\infty<\alpha$ and
 $\max\{\alpha,1/2\}<\beta\leq1,$ $Q_{\alpha}^{\beta}(\mathbb{R}^{n})$ is decreasing in $\alpha$ for a fixed $\beta.$ Moreover, if $\alpha\in(-\infty,\beta-1),$ then all
 $Q_{\alpha}^{\beta}(\mathbb{R}^{n})$ equal to $Q^{\beta}_{-\frac{n}{2}+\beta-1}(\mathbb{R}^{n}):=BMO^{\beta}(\mathbb{R}^{n}).$
 Thus, when $k=0$ and $\alpha=-\frac{n}{2}+\beta-1,$  (\ref{mean ch Q a b 14}) implies a special JN type inequality, that is,
 for $f\in L^{2}(\mathbb{R}^{n})\cap BMO^{\beta}(\mathbb{R}^{n})$ and $t\leq\|f\|_{BMO^{\beta}(\mathbb{R}^{n})},$
\begin{equation} \label{mean ch Q a b 26}
|\{x\in \mathbb{R}^{n}: |f|>t\}|\leq
\frac{B\|f\|_{L^{2}(\mathbb{R}^{n})}^{2}}{t^{2}}
\exp\left(\frac{-bt}{\|f\|_{BMO^{\beta}(\mathbb{R}^{n})}}\right).
\end{equation}
 When $t>\|f\|_{BMO^{\beta}(\mathbb{R}^{n})},$ we get a weaker form  of (\ref{mean ch Q a b
 26}).

\begin{proposition}\label{mean ch Q a b theorem4} Let  $\beta\in (1/2,1].$
 If $f\in BMO^{\beta}(\mathbb{R}^{n})\cap L^{2}(\mathbb{R}^{n}),$ then

  \item{(i)}  (\ref{mean ch Q a b 26}) holds
 for all  $t\leq\|f\|_{BMO^{\beta}(\mathbb{R}^{n})};$

\item{(ii)} \begin{equation} \label{mean ch Q a b 201'}
|\{x\in\mathbb{R}^{n}: f(x)>t\}|\leq
\frac{B\|f\|_{L^{2}(\mathbb{R}^{2})}^{2}}{\|f\|^{2}_{BMO^{\beta}(\mathbb{R}^{n})}}
\end{equation}
 holds for  all $t>\|f\|_{BMO^{\beta}(\mathbb{R}^{n})}.$
\end{proposition}

 When $\beta=1$  and  $t>\|f\|_{BMO(\mathbb{R}^{n})},$ (\ref{mean ch Q a b 26}) also holds  and implies the following GN type inequalities in $Q_{\alpha}(\mathbb{R}^{n})$
 which can also be deduced from  \cite[Theorem 2]{Jiecheng Chen Xiangrong Zhu} and
\cite[Theorem 2.3]{M Essen S Janson L Peng J Xiao}: for
$-\infty<\alpha<1$ and $1\leq r\leq p<\infty,$
\begin{equation}\label{mean ch Q a b 21}
\|f\|_{L^{p}(\mathbb{R}^{n})}\leq
C_{n}p\|f\|_{L^{r}(\mathbb{R}^{n})}^{r/p}\|f\|_{Q_{\alpha}(\mathbb{R}^{n})}^{1-r/p},
\end{equation}
for   $f\in L^{r}(\mathbb{R}^{n})\cap Q_{\alpha}(\mathbb{R}^{n}).$
Here,  $C_{*,\cdots,*}$ denotes a constant which depends  only on
the quantities appearing in the subscript indexes.

  As an application of (\ref{mean ch Q a b 21}), we establish the Trudinger-Moser type inequality which implies a
 generalized JN type inequality.

\begin{theorem}\label{mean ch Q a b corollary1}

\item{(i)}   There exists a positive constant $\gamma_{n}$  such that  for every $0<\zeta<\gamma_{n}$
\begin{equation}\label{mean ch Q a b 22}
\int_{\mathbb{R}^{n}}\Phi_{p}\left(\zeta
\left(\frac{|f(x)|}{\|f\|_{Q_{\alpha}(\mathbb{R}^{n})}}\right)\right)dx
\leq
C_{n,\zeta}\left(\frac{\|f\|_{L^{p}(\mathbb{R}^{n})}}{\|f\|_{Q_{\alpha}(\mathbb{R}^{n})}}\right)^{p}
\end{equation}
 holds for all
 $$f\in L^{p}(\mathbb{R}^{n})\cap Q_{\alpha}(\mathbb{R}^{n})\quad\hbox{with}\quad  1<  p<\infty \quad\hbox{and}\quad   -\infty<\alpha<1.$$
  Here $\Phi_{p}$ is the function defined by
$$
\Phi_{p}(t)= e^{t}-\sum_{j<p,j\in
\mathbb{N}\cup\{0\}}\frac{t^{j}}{j!}, t\in \mathbb{R}.
$$

\item{(ii)} There exists a  positive constant $\gamma_{n}$  such that
\begin{equation}\label{mean ch Q a b 23}
|\{x\in \mathbb{R}^{n}: |f|>t\}|\leq C_{n}
\frac{\|f\|_{L^{2}(\mathbb{R}^{n})}^{2}}{\|f\|_{Q_{\alpha}(\mathbb{R}^{n})}^{2}}
\frac{1}{\left(\exp\left(\frac{t\gamma_{n}}{\|f\|_{Q_{\alpha}(\mathbb{R}^{n})}}\right)
-1-\frac{t\gamma_{n} }{\|f\|_{Q_{\alpha}(\mathbb{R}^{n})}}\right)}
\end{equation}
 holds for all $t>0$ and
 $$f\in L^{2}(\mathbb{R}^{n})\cap Q_{\alpha}(\mathbb{R}^{n})\quad\hbox{with}
\quad -\infty<\alpha<1.$$
 In particular, we have
\begin{equation}\label{mean ch Q a b 24}
|\{x\in \mathbb{R}^{n}: |f|>t\}|\leq C_{n}
\frac{\|f\|_{L^{2}(\mathbb{R}^{n})}^{2}}{\|f\|_{Q_{\alpha}(\mathbb{R}^{n})}^{2}}
{\exp\left(-\frac{t\gamma_{n}}{\|f\|_{Q_{\alpha}(\mathbb{R}^{n})}}\right)}
\end{equation}
 holds for all $t>\|f\|_{Q_{\alpha}(\mathbb{R}^{n})}$ and
 $$f\in L^{2}(\mathbb{R}^{n})\cap
 Q_{\alpha}(\mathbb{R}^{n})\quad
\hbox{with}\quad-\infty<\alpha<1.$$
\end{theorem}

 We can also get the following Brezis-Gallouet-Wainger type inequalities.

\begin{proposition}\label{mean ch Q a b theorem6}
 For every  $1<q<\infty$ and $n/q<s<\infty,$ we have
\begin{equation}\label{mean ch Q a b 25}
\|f\|_{L^{\infty}(\mathbb{R}^{n})}\leq
C_{n,p,q,s}\left(1+(\|f\|_{L^{p}(\mathbb{R}^{n})}+\|f\|_{Q_{\alpha}(\mathbb{R}^{n})})
\log(e+\|(-\triangle)^{s/2}f\|_{L^{q}(\mathbb{R}^{n})})\right)
\end{equation}
 holds for all
 $(-\triangle)^{s/2}f\in L^{q}(\mathbb{R}^{n})$
 satisfying
 $$
 f\in L^{p}(\mathbb{R}^{n})\cap Q_{\alpha}(\mathbb{R}^{n})
 \quad \hbox{when}\quad 1\leq p<\infty\quad\hbox{and}\quad -\infty<\alpha<1.$$
  \end{proposition}

 In the next section, we prove our main results. We verify Propositions \ref{mean ch Q a b theorem2}-\ref{mean ch Q a b theorem3} for $\beta\in(1/2,1]$  by applying similar arguments in
  the proof of Yue and Dafni \cite[Theorems 1-2]{H Yue G Dafni} for $\beta=1.$ We deduce Proposition \ref{mean ch Q a b theorem4} from a special case of
  Proposition \ref{mean ch Q a b theorem3}. Finally, we demonstrate  Theorem \ref{mean ch Q a b corollary1} and Proposition \ref{mean ch Q a b theorem6}
   by applying (\ref{mean ch Q a b 21}) and the $L^{p}-L^{q}$ estimates for $e^{-t(-\triangle)^{s/2}}.$

\section{Proofs of Main Results}

\subsection{Proof of Proposition \ref{mean ch Q a b theorem2}}

According to Proposition \ref{mean ch Q a b theorem1}, it suffices
to prove that $\Psi_{f,\alpha,\beta}(I)$ is bounded independent of
$I.$  More specially, we will prove  for any $p<q,$ we have
\begin{equation}
\Psi_{f,\alpha,\beta}^{q}(I):=(l(I))^{4\beta-4}
\sum_{k=0}^{\infty}2^{(2(\alpha-\beta+1)-n)k}\sum_{J\in
\mathcal{D}_{k}(I)}\Phi_{f}^{q}(J)\leq BK_{C,c,q,p},
\end{equation}
where $B,C,c$ are the constants appearing in (\ref{mean ch Q a b
113}), and $K_{C,c,q,p}$ is a constant depending only on $C,c,p,$
and $q.$ When $q=2,$
$\Psi_{f,\alpha,\beta}^{q}(I)=\Psi_{f,\alpha,\beta}(I),$ so this
implies the theorem.

 For a fixed cube  $I,$ and any $J\in \mathcal{D}_{k}(I),$
 let $\int_{J}|f(x)-f(J)|^{q}dx=q\int_{0}^{\infty}t^{q-1}m_{J}(t)dt.$
Using the Monotone Convergence Theorem and the inequality (\ref{mean
ch Q a b 113}), we have
\begin{eqnarray*}
\Psi_{f,\alpha,\beta}^{q}(I)
&=&(l(I))^{4\beta-4}\sum_{k=0}^{\infty}2^{(2(\alpha-\beta+1)-n)k}\sum_{J\in
\mathcal{D}_{k}(t)}\frac{q}{|J|}\int_{0}^{\infty}t^{q-1}m_{J}(t)dt\\
&=&q\int_{0}^{\infty}t^{q-1}\left((l(I))^{4\beta-4}
\sum_{k=0}^{\infty}2^{(2(\alpha-\beta+1)-n)k}\sum_{J\in
\mathcal{D}_{k}(I)}\frac{m_{J}(t)}{|J|}\right)dt\\
&\leq&
q\int_{0}^{\infty}t^{q-1}B(1+\left(\frac{C}{t}\right)^{p})e^{-ct}dt\\
&=&qB\left(c^{-q}\int_{0}^{\infty}u^{q-1}e^{-u}du+C^{p}c^{-(q-p)}\int_{0}^{\infty}u^{q-p-1}e^{-u}du\right)\\
&=&qB(c^{-q}\Gamma(q)+C^{p}c^{-(q-p)}\Gamma(q-p))
\end{eqnarray*}
where $\Gamma(y)=\int_{0}^{\infty}u^{y-1}e^{-u}du.$ Since $0\leq
p<q,$ $\Gamma(q)$ and $\Gamma(q-p)$ are finite. Thus, we can get the
desired inequality by taking
$K_{C,c,p,q}=q(c^{-q}\Gamma(q)+C^{p}c^{-(q-p)}\Gamma(q-p)).$

\subsection{Proof of Proposition \ref{mean ch Q a b theorem3}}
Assume that $f$ is a nontrivial element of
$Q_{\alpha}^{\beta}(\mathbb{R}^{n}).$ Then
$\gamma=\sup_{I}(\Psi_{f,\alpha,\beta}(I))^{1/2}<\infty.$ For all
cubes $I$ we have
\begin{equation}\label{mean ch Q a b 15}
(l(I))^{2\beta-2}\frac{1}{|I|}\int_{I}|f(x)-f(I)|dx\leq
((l(I))^{4\beta-4}\Phi_{f}^{2}(I))^{1/2}\leq
(\Psi_{f,\alpha,\beta}(I))^{1/2}\leq \gamma.
\end{equation}
For a cube $I$ and each $J\in\mathcal{D}_{k}(I),$ we have by the
Chebyshev inequality, for $t>0,$
$$m_{J}(t)\leq t^{-2}\int_{J}|f(x)-f(J)|^{2}dx.$$
Thus we get
\begin{equation}\label{mean ch Q a b 16}
(l(I))^{4\beta-4}\sum_{k=0}^{\infty}2^{(2(\alpha-\beta+1)-n)k}\sum_{J\in
\mathcal{D}_{k}(I)}\frac{m_{J}(t)}{|J|}\leq
t^{-2}\Psi_{f,\alpha,\beta}(I)\leq t^{-2}\gamma^{2}.
\end{equation}
Thus, if $t\leq \gamma,$ then (\ref{mean ch Q a b 14}) holds with
$B=e$ and $b=1.$

To consider the case of $t>\gamma,$ we need the Calder\'{o}n-Zygmund
decomposition, see Calder\'{o}n and Zygmund \cite{Calderon Zygmund},
 and Neri \cite{Neri}.

\begin{lemma}\label{mean ch Q a b lemma6}
Assume that $f$ is a nonnegative function in $L^{1}(\mathbb{R}^{n})$
and $\xi$ is a positive constant. There is a decomposition
$\mathbb{R}^{n}=P\cup \Omega,$ $P\cap\Omega=\emptyset,$ such that\\
(a) $\Omega=\cup_{k=1}^{\infty}I_{k},$ where $I_{k}$ is a collection
of cubes whose interiors are disjoint;\\
(b) $f(x)\leq \xi$ for a.e. $x\in P;$\\
 (c) $\xi<\frac{1}{|I|}\int_{I}f(x)dx\leq 2^{n}\xi,$ for all $I$ in
 the collection $\{I_{k}\}.$\\
 (d) $\xi|\triangle|\leq \int_{\triangle}f(x)dx\leq
 2^{n}\xi|\triangle|,$ if $\triangle$ is any union of cubes $I$
 from $\{I_{k}\}.$
  \end{lemma}

In the following we  fix a cube  $I.$  For $\xi=t(l(I))^{2-2\beta}$
with any $t>0,$ we apply the Calder\'{o}n-Zygmund decomposition to
$|f(x)-f(J)|$  on  a subcube $J\in \mathcal{D}_{k}(I).$  Set
$\Omega=\Omega_{J}(t),$ $P=J\backslash\Omega_{J}(t).$

From Cauchy-Schwarz inequality and (d) of Lemma \ref{mean ch Q a b
lemma6},  we get
\begin{equation}\label{varant of d}
(t(l(I))^{2-2\beta})^{2}|\triangle|\leq
\int_{\triangle}|f(x)-f(J)|^{2}dx
\end{equation}
 for any union $\triangle$ of the cubes $K$ in the
decomposition of $\Omega_{J}(t).$ Inequality  (\ref{varant of d})
with $\triangle=\Omega_{J}(t)$ gives us  a variant  of inequality
(\ref{mean ch Q a b 16}):
\begin{equation}\label{mean ch Q a b 116}
(l(I))^{4\beta-4}\sum_{k=0}^{\infty}2^{(2(\alpha-\beta+1)-n)k}\sum_{J\in
\mathcal{D}_{k}(I)}\frac{|\Omega_{J}(t)|}{|J|}\leq
\frac{\Psi_{f,\alpha,\beta}(I)}{(t(l(I))^{2-2\beta})^{2}}\leq
\left(\frac{\gamma}{(t(l(I))^{2-2\beta})}\right)^{2}
\end{equation}
for all $t>0.$

When  $t\geq\gamma,$ we can strengthen  the estimate (c) in Lemma
\ref{mean ch Q a b lemma6}  as follows:
\begin{equation}\label{strengthened}
t(l(I))^{2-2\beta}<\frac{1}{|K|}\int_{K}|f(x)-f(J)|dx\leq
 (2^{n}\gamma+t)(l(I))^{2-2\beta}
  \end{equation}
for all cubes $K$ in the decomposition of $\Omega_{J}(t).$ In fact,
note that  $K$ is such a cube, then $K\neq J.$ Otherwise, (\ref{mean
ch Q a b 15}) implies
$$\frac{1}{|J|}\int_{J}|f(x)-f(J)|dx
\leq \gamma(l(I))^{2-2\beta}\leq t(l(I))^{2-2\beta}.$$ This
contradicts (c). It follows from the proof of the
Calder\'{o}n-Zygmund decomposition (see, Stein \cite{Stein} ) that
$K$ must have a ``parent" cube $K^{*}\subset J$ satisfying
$K\in\mathcal{D}_{1}(K^{*}),$ $l(K^{*})=2l(K)$ and
$$|f(K^{*})-f(J)|\leq |K^{*}|^{-1}\int_{K^{*}}|f(x)-f(J)|dx\leq t(l(I))^{2-2\beta}.$$
Then (\ref{mean ch Q a b 15}) implies
\begin{eqnarray*}
t(l(I))^{2-2\beta}<\frac{1}{|K|}\int_{K}|f(x)-f(J)|dx
&\leq&\frac{1}{|K|}\int_{K}|f(x)-f(K^{*})|dx+|f(K^{*})-f(J)|\\
&\leq&\frac{2^{n}}{|K^{*}|}\int_{K^{*}}|f(x)-f(K^{*})|dx+t(l(I))^{2-2\beta}\\
&\leq&(2^{n}\gamma+t)(l(I))^{2-2\beta}.
\end{eqnarray*}

 There holds $\Omega_{J}(t')\subset\Omega_{J}(t)$ for
$0<t<t'.$ In fact, for any cube $K\in
\Omega_{J}(t')\backslash\Omega_{J}(t),$ we get $K\subset
J\backslash\Omega_{J}(t).$ So,  property (b) tells us
$$t(l(I))^{2-2\beta}\geq \frac{1}{|K|}\int_{K}|f(x)-f(J)|dx>t'(l(I))^{2-2\beta}.$$
This  is a contradiction.

Letting $t'=t+2^{n+1}\gamma$ for $t\geq \gamma,$ we claim that
\begin{equation}\label{impoartant rewrite}
|\Omega_{J}(t')|\leq2^{-n}|\Omega_{J}(t)|.
\end{equation}To prove this, take a cube $K$ in the decomposition for
$\Omega_{J}(t).$ Then (\ref{strengthened}) implies that
$$\frac{1}{|K|}\int_{K}|f(x)-f(J)|dx\leq(2^{n}\gamma+t)
(l(I))^{2-2\beta}<t'(l(I))^{2-2\beta}.$$ Thus, $K$ is not  a cube in
the decomposition of $\Omega_{J}(t'),$
 and was further subdivided. Set
$\triangle'=K\cap\Omega_{J}(t').$ If $\triangle'\neq\emptyset,$ it
must be a union of cubes from the decomposition of $\Omega_{J}(t').$
Thus, according to (d) of Lemma \ref{mean ch Q a b lemma6},
(\ref{mean ch Q a b 15}) and (\ref{strengthened}),
\begin{eqnarray*}
t'(l(I))^{2-2\beta}&\leq& |\triangle'|^{-1}\int_{\triangle'}|f(x)-f(J)|dx\\
&\leq&|\triangle'|^{-1}\int_{\triangle'}|f(x)-f(K)|dx+|f(K)-f(J)|\\
&\leq&|\triangle'|^{-1}|K|
\frac{1}{|K|}\int_{\triangle'}|f(x)-f(K)|dx
+\frac{1}{|K|}\int_{K}|f(x)-f(J)|dx\\
&\leq&|\triangle'|^{-1}|K|\gamma(l(K))^{2-2\beta}
+(2^{n}\gamma+t)(l(I))^{2-2\beta}\\
&\leq&|\triangle'|^{-1}|K|\gamma(l(I))^{2-2\beta}
+(2^{n}\gamma+t)(l(I))^{2-2\beta}
\end{eqnarray*}
since $2-2\beta>0$ and $K\subset I.$ Replacing $t'$ by
$t+2^{n+1}\gamma,$ dividing by $(l(I))^{2-2\beta},$ subtracting $t$
and dividing by $\gamma,$ we have
$$
(2^{n+1}-2^{n})\leq |\triangle'|^{-1}|K| \quad\hbox{and}\quad |K\cap
\Omega_{J}(t')|=|\triangle'|\leq 2^{-n}|K|
$$
 for any cube $K$ in the decomposition of $\Omega_{J}(t).$
Summing over all such $K,$ and noting  that
$\Omega_{J}(t')=\Omega_{J}(t)\cap\Omega_{J}(t'),$ we prove
 (\ref{impoartant rewrite}).

 For each $J\in \mathcal{D}_{k}(I),$ property (b) of the
 decomposition for $|f-f(J)|$ implies that
 \begin{equation}\label{mean ch Q a b 18}
m_{J}(t(l(I))^{2-2\beta})=|\{x\in J:
|f(x)-f(J)|>t(l(I))^{2-2\beta}\}|\leq |\Omega_{J}(t)|.
 \end{equation}

For $t>\gamma,$ let $j$ be the integer part of
$\frac{t-\gamma}{2^{n+1}\gamma}$ and $s=(1+j2^{n+1})\gamma.$ Then
$\gamma\leq s\leq t.$ Thus  one obtains from (\ref{mean ch Q a b
18}) that
\begin{eqnarray*}
&&(l(I))^{4\beta-4}
\sum_{k=0}^{\infty}2^{(2(\alpha-\beta+1)-n)k}\sum_{J\in\mathcal{D}_{k}(I)}\frac{m_{J}(t)}{|J|}\\
&=&(l(I))^{4\beta-4}
\sum_{k=0}^{\infty}2^{(2(\alpha-\beta+1)-n)k}\sum_{J\in\mathcal{D}_{k}(I)}
\frac{m_{J}((l(I))^{2-2\beta}t(l(I))^{2\beta-2})}{|J|}\\
&\leq&(l(I))^{4\beta-4}
\sum_{k=0}^{\infty}2^{(2(\alpha-\beta+1)-n)k}
\sum_{J\in\mathcal{D}_{k}(I)}\frac{m_{J}((l(I))^{2-2\beta}s(l(I))^{2\beta-2})}{|J|}\\
&\leq&(l(I))^{4\beta-4}\sum_{k=0}^{\infty}2^{(2(\alpha-\beta+1)-n)k}
\sum_{J\in\mathcal{D}_{k}(I)}\frac{|\Omega_{J}((1+j2^{n+1})\gamma(l(I))^{2\beta-2})|}{|J|}\\
&\leq&(l(I))^{4\beta-4}\sum_{k=0}^{\infty}2^{(2(\alpha-\beta+1)-n)k}
\sum_{J\in\mathcal{D}_{k}(I)}\frac{|\Omega_{J}(\gamma(l(I))^{2\beta-2}+j2^{n+1}\gamma)|}{|J|}\\
&\leq&2^{-n}(l(I))^{4\beta-4}\sum_{k=0}^{\infty}2^{(2(\alpha-\beta+1)-n)k}
\sum_{J\in\mathcal{D}_{k}(I)}\frac{|\Omega_{J}(\gamma(l(I))^{2\beta-2}+(j-1)2^{n+1}\gamma)|}{|J|}
\end{eqnarray*}
 if
$(l(I))^{2\beta-2}\geq 1,$ by  using
 (\ref{impoartant rewrite}) for
$$
t=((l(I))^{2\beta-2}+(j-1)2^{n+1})\gamma\quad \hbox{and} \quad
t'=((l(I))^{2\beta-2}+j2^{n+1})\gamma.
$$
 Iterating the previous estimate
$j$ times and using (\ref{mean ch Q a b 116}) with
$t=\gamma(l(I))^{2\beta-2},$ one has
\begin{eqnarray*}
&&(l(I))^{4\beta-4}\sum_{k=0}^{\infty}2^{(2(\alpha-\beta+1)-n)k}\sum_{J\in\mathcal{D}_{k}(I)}
\frac{m_{J}(t)}{|J|}\\
&\leq
&2^{-nj}(l(I))^{4\beta-4}\sum_{k=0}^{\infty}2^{(2(\alpha-\beta+1)-n)k}\sum_{J\in\mathcal{D}_{k}(I)}
\frac{|\Omega_{J}(\gamma(l(I))^{2\beta-2})|}{|J|}\\
&\leq&2^{-nj}\gamma^{2}\gamma^{-2}\\
&\leq& 2^{-n\left(\frac{t-\gamma}{2^{n+1}\gamma}-1\right)}\\
&=&2^{-\frac{n}{2^{n+1}}(t/\gamma)}2^{\frac{n}{2^{n+1}}+n}.
\end{eqnarray*}
Taking $B=2^{n/2^{n+1}+n}$ and $b=\frac{n}{2^{n+1}}\ln 2,$ we get
(\ref{mean ch Q a b 14}) when $(l(I))^{2\beta-2}\geq 1.$

If  $(l(I))^{2\beta-2}< 1,$ using (\ref{mean ch Q a b 18}) and
(\ref{varant of d}), one has
\begin{eqnarray*}
&&(l(I))^{4\beta-4}
\sum_{k=0}^{\infty}2^{(2(\alpha-\beta+1)-n)k}\sum_{J\in\mathcal{D}_{k}(I)}\frac{m_{J}(t)}{|J|}\\
&\leq&(l(I))^{4\beta-4}\sum_{k=0}^{\infty}2^{(2(\alpha-\beta+1)-n)k}
\sum_{J\in\mathcal{D}_{k}(I)}\frac{|\Omega_{J}(t(l(I))^{2\beta-2})|}{|J|}\\
&\leq&\gamma^{2}t^{-2} \leq1
\end{eqnarray*}
which yields  (\ref{mean ch Q a b 14''}).

\subsection{Proof of Proposition \ref{mean ch Q a b theorem4}}
Taking $k=0$  and $\alpha=-\frac{n}{2}+\beta-1$ in  (\ref{mean ch Q
a b 14}), we get that
$$
(l(I))^{4\beta-4}\frac{m_{I}(t)}{|I|}\leq
B\frac{\|f\|^{2}_{BMO^{\beta}(\mathbb{R}^{n})}}{t^{2}}\exp\left(\frac{-bt}{\|f\|_{BMO^{\beta}(\mathbb{R}^{n})}}\right)
$$
holds for $t\leq\|f\|_{BMO^{\beta}(\mathbb{R}^{n})}$ and any cube
$I.$ Thus for $t\leq\|f\|_{BMO^{\beta}(\mathbb{R}^{n})}$ and any
cube $I,$ we have
\begin{eqnarray*}
&&(l(I))^{4\beta-4}\frac{m_{I}(t)}{|I|}\int_{I}|f(x)-f(I)|^{2}dx\\&\leq&
B\frac{\|f\|^{2}_{BMO^{\beta}(\mathbb{R}^{n})}}{t^{2}}\exp\left(\frac{-bt}{\|f\|_{BMO^{\beta}(\mathbb{R}^{n})}}\right)\int_{I}|f(x)-f(I)|^{2}dx\\
&\leq&B\frac{\|f\|^{2}_{BMO^{\beta}(\mathbb{R}^{n})}}{t^{2}}\exp\left(\frac{-bt}{\|f\|_{BMO^{\beta}(\mathbb{R}^{n})}}\right)\int_{I}|f(x)|^{2}dx
\\
&\leq&B\frac{\|f\|^{2}_{BMO^{\beta}(\mathbb{R}^{n})}}{t^{2}}\exp\left(\frac{-bt}{\|f\|_{BMO^{\beta}(\mathbb{R}^{n})}}
\right)\int_{\mathbb{R}^{n}}|f(x)|^{2}dx.
\end{eqnarray*}
This tells us
\begin{eqnarray}\label{haha}
&&m_{I}(t)
\frac{(l(I))^{4\beta-4}}{|I|}\int_{I}|f(x)-f(I)|^{2}dx\\
&\leq&
B\frac{\|f\|^{2}_{BMO^{\beta}(\mathbb{R}^{n})}}{t^{2}}\exp\left(\frac{-bt}{\|f\|_{BMO^{\beta}(\mathbb{R}^{n})}}
\right)\int_{\mathbb{R}^{n}}|f(x)|^{2}dx.
\end{eqnarray}
 According to the
definition of $BMO^{\beta}(\mathbb{R}^{n}),$ see Li and Zhai
\cite{Pengtao Li Zhicun Zhai}, we have
$$f\in BMO^{\beta}(\mathbb{R}^{n})\Longleftrightarrow
\|f\|^{2}_{BMO^{\beta}(\mathbb{R}^{n})}=
\sup_{I}\frac{(l(I))^{4\beta-4}}{|I|}\int_{I}|f(x)-f(I)|^{2}dx<\infty.$$
Thus,  we get
\begin{eqnarray*}
&&m_{I}(t)\|f\|^{2}_{BMO^{\beta}(\mathbb{R}^{n})}\\
&\leq& B\frac{\|f\|^{2}_{BMO^{\beta}(\mathbb{R}^{n})}}{t^{2}}
\exp\left(\frac{-bt}{\|f\|_{BMO^{\beta}(\mathbb{R}^{n})}}\right)
\int_{\mathbb{R}^{n}}|f(x)|^{2}dx,
\end{eqnarray*}
 for $t\leq\|f\|_{BMO^{\beta}(\mathbb{R}^{n})}.$
Then, taking an increasing sequence of cubes covering
$\mathbb{R}^{n},$  we obtain
\begin{equation} \label{mean ch Q a b 20}
|\{x\in\mathbb{R}^{n}: |f(x)|>t\}|\leq \frac{B}{t^{2}}
\exp\left(\frac{-bt}{\|f\|_{BMO^{\beta}(\mathbb{R}^{n})}}\right)
\int_{\mathbb{R}^{n}}|f(x)|^{2}dx
\end{equation}
 for $t\leq\|f\|_{BMO^{\beta}(\mathbb{R}^{n})},$
since $f(I)\longrightarrow0$ as $l(I)\longrightarrow\infty.$
 Finally, we get (\ref{mean ch Q a b 26}). Similarly, we can prove  (\ref{mean ch Q a b 201'}) since
$\exp\left(\frac{-bt}{\|f\|_{BMO^{\beta}(\mathbb{R}^{n})}}\right)\leq
1$  for $t>\|f\|_{BMO^{\beta}(\mathbb{R}^{n})}.$

\subsection{Proof of Proposition \ref{mean ch Q a b corollary1}}
(i) According to (\ref{mean ch Q a b 21}), we have
\begin{eqnarray*}
\int_{\mathbb{R}^{n}}\Phi_{p,r}
\left(\zeta\frac{|f(x)|}{\|f\|_{Q_{\alpha}(\mathbb{R}^{n})}}\right)dx&=&
\int_{\mathbb{R}^{n}}\sum_{j\geq p,j\in
\mathbb{N}}\frac{\zeta^{j}}{j!}
\left(\frac{|f(x)|}{\|f\|_{Q_{\alpha}(\mathbb{R}^{n})}}\right)^{j}dx\\
&\leq&\sum_{j\geq p,j\in \mathbb{N}}\frac{\zeta^{j}}{j!}
\frac{\|f\|_{L^{j}(\mathbb{R}^{n})}^{j}}{\|f\|_{Q_{\alpha}(\mathbb{R}^{n})}^{j}}
 \\ &\leq&\sum_{j\geq p,j\in
\mathbb{N}}\frac{\zeta^{j}}{j!}
\frac{\left(C_{n}j\|f\|_{L^{p}(\mathbb{R}^{n})}^{p/j}\|f\|_{Q_{\alpha}(\mathbb{R}^{n})}^{1-p/j}\right)^{j}}
{\|f\|_{Q_{\alpha}(\mathbb{R}^{n})}^{j}}\\
&\leq& \sum_{j\geq p,j\in \mathbb{N}}a_{j}(\zeta C_{n})^{j}
\left(\frac{\|f\|_{L^{p}(\mathbb{R}^{n})}}{\|f\|_{Q_{\alpha}(\mathbb{R}^{n})}}\right)^{p}
\end{eqnarray*}
with $a_{j}=\frac{j^{j}}{j!}.$ Since
$\lim_{j\longrightarrow\infty}\frac{a_{j}}{a_{j+1}}=e^{-1},$ the
power series of the above right hand side converges provided $\zeta
C_{n}<e^{-1}$ i.e. $\zeta<\gamma_{n}:=(C_{n}e)^{-1}.$\\
(ii) According to (i) with $p=2,$ we have
$$
\int_{\mathbb{R}^{n}}\left(\exp\left(\gamma_{n}\frac{|f(x)|}{\|f\|_{Q_{\alpha}(\mathbb{R}^{n})}}
\right)-1-\gamma_{n}\frac{|f(x)|}{\|f\|_{Q_{\alpha}(\mathbb{R}^{n})}}\right)dx
\leq
C_{n}\frac{\|f\|_{L^{2}(\mathbb{R}^{n})}^{2}}{\|f\|_{Q_{\alpha}(\mathbb{R}^{n})}^{2}}
.$$  On the other hand, since the distribution function
$m(t)=|\{x\in \mathbb{R}^{n}: |f(x)|>t\}|$ is non-increasing, we
have
\begin{eqnarray*}
&&\int_{\mathbb{R}^{n}}\left(\exp\left(\gamma_{n}\frac{|f(x)|}
{\|f\|_{Q_{\alpha}(\mathbb{R}^{n})}}\right)-1-\gamma_{n}\frac{|f(x)|}
{\|f\|_{Q_{\alpha}(\mathbb{R}^{n})}}\right)dx\\
&=&\sum_{j=2}^{\infty}\frac{\gamma_{n}^{j}}{j!}\frac{\|f\|_{L^{j}(\mathbb{R}^{n})}^{j}}{\|f\|^{j}_{Q_{\alpha}(\mathbb{R}^{n})}}\\
&=&\sum_{j=2}^{\infty}\frac{\gamma_{n}^{j}}{j!}\frac{j}{\|f\|^{j}_{Q_{\alpha}(\mathbb{R}^{n})}}\int_{0}^{\infty}m(s)s^{j-1}ds\\
&\geq&m(t)\sum_{j=2}^{\infty}\frac{\gamma_{n}^{j}}{j!}\frac{j}{\|f\|^{j}_{Q_{\alpha}(\mathbb{R}^{n})}}\int_{0}^{t}s^{j-1}ds\\
&=&m(t)\sum_{j=2}^{\infty}\frac{1}{j!}\left(\frac{\gamma_{n} t}
{\|f\|_{Q_{\alpha}(\mathbb{R}^{n})}}\right)^{j}\\
&=&m(t)\left(\exp\left(\frac{\gamma_{n}
t}{\|f\|^{j}_{Q_{\alpha}(\mathbb{R}^{n})}}\right)-1-\frac{\gamma_{n}
t} {\|f\|_{Q_{\alpha}(\mathbb{R}^{n})}}\right)
\end{eqnarray*}
for all $t>0.$ Thus, we have
$$m(t)\leq C_{n}\frac{\|f\|_{L^{2}(\mathbb{R}^{n})}^{2}}{\|f\|_{Q_{\alpha}(\mathbb{R}^{n})}^{2}}
\frac{1}{ \left(\exp\left(\frac{\gamma_{n}
t}{\|f\|^{j}_{Q_{\alpha}(\mathbb{R}^{n})}}\right)-1-\frac{\gamma_{n}
t} {\|f\|_{Q_{\alpha}(\mathbb{R}^{n})}}\right)}.$$

\subsection{Proof of Proposition \ref{mean ch Q a b theorem6}}
  We will use some facts about the factional heat equations  $$
 \partial_{t}v(t,x)+(-\triangle)^{s/2}v(t,x)=0\quad \hbox{for}\quad (t,x)\in (0,\infty)\times \mathbb{R}^{n} $$
 with initial data $v(0,x)=g(x)$ for $x\in \mathbb{R}^{n}.$ Here  $\mathcal{F}((-\triangle)^{s/2}v(t,x))(\xi)=|\xi|^{s}\mathcal{F}v(t,\xi)$ and
  $v_{g}(t,x)=e^{-(\triangle)^{s/2}}g(x)=K_{t}^{s}(x)\ast g(x)$
 with $K_{t}^{s}(\cdot)=\mathcal {F}^{-1}(e^{-t|\cdot|^{s}})$
  where $\mathcal{F}$ and  $\mathcal{F}^{-1}$ denote the Fourier  transformation and its
  inverse. We need the $L^{p}\longrightarrow L^{q}$ estimates for
  the semigroup $\{e^{-t(-\triangle)^{s/2}}\}_{t\geq 0}.$ For the proof,
   see, for example,  Kozono-Wadade \cite[Lemma 3.4]{H Kozono H
Wadade} or Miao-Yuan-Zhang \cite[Lemma 3.1]{C. Miao}.

\begin{lemma}\label{mean ch Q a b lemma7}
For every  $0<s<\infty ,$  there exists a constant $C_{n,s}$
depending only on $n$ and $s$ such that
$$\|e^{-t(-\triangle)^{s/2}}g\|_{L^{q}(\mathbb{R}^{n})}\leq C_{n,s}t^{-\frac{n}{s}
\left(\frac{1}{p}-\frac{1}{q_{1}}\right)}\|g\|_{L^{p}(\mathbb{R}^{n})}.$$
 holds for all $g\in L^{p}(\mathbb{R}^{n}),$  $t>0$ and $1\leq p\leq q\leq \infty.$
  \end{lemma}

For any $g(x)$ in the  Schwartz class of rapidly
   decreasing functions   $\mathscr{S}(\mathbb{R}^{n}),$ define
   $v_{g}(t,x)=e^{-(\triangle)^{s/2}}g(x)$
as the solution of fractional heat equation
   $$
 \partial_{t}v(t,x)+(-\triangle)^{s/2}v(t,x)=0 $$
 with initial data $g.$
Fix $f\in L^{2}(\mathbb{R}^{n})\cap
Q_{\alpha}^{\beta}(\mathbb{R}^{n})$ with $(-\triangle)^{s/2}f\in
L^{q}.$ Then
\begin{eqnarray*}
\int_{0}^{t}\langle-(-\triangle)^{s/2}f(x), v(s,x)\rangle
ds&=&\int_{0}^{t}\langle f(x), -(-\triangle)^{s/2}v(s,x)\rangle ds\\
&=&\int_{0}^{t}\langle f(x), \partial_{s} v(s,x)\rangle dt\\
&=& \langle f(x),v(t,x)\rangle-\langle f(x),g(x)\rangle.
\end{eqnarray*}
Thus
$$|\langle f, g\rangle|\leq |\langle f(x),
v(t,x)\rangle|+\int_{0}^{t}|\langle (-\triangle)^{s/2}f(x),
v(s,x)\rangle|ds=I_{1}+I_{2}
$$
for all $t>0.$ Here $\langle\cdot, \cdot\rangle$ denote the
inner-product in $L^{2}.$
 Thus
H\"{o}lder inequality,  Lemma \ref{mean ch Q a b lemma7} and
(\ref{mean ch Q a b 21}) imply that
\begin{eqnarray*}
I_{1}&\leq&
\|f\|_{L^{q_{1}(\mathbb{R}^{n})}}\|v(t,\cdot)\|_{L^{q'_{1}}(\mathbb{R}^{n})}=\|f\|_{L^{q_{1}}(\mathbb{R}^{n})}
\|e^{-t(-\triangle)^{s/2}}g\|_{L^{q'_{1}}(\mathbb{R}^{n})}\\
&\leq&C_{n,s}q_{1}t^{-\frac{n}{sq_{1}}}(\|f\|_{L^{p}(\mathbb{R}^{n})}+\|f\|_{Q_{\alpha}^{\beta}(\mathbb{R}^{n})})
\|g\|_{L^{1}(\mathbb{R}^{n})}
\end{eqnarray*}
for all $t>0$ and $p\leq q_{1}< \infty.$ Similarly, we have
\begin{eqnarray*}
I_{2}&\leq&\int^{t}_{0}\|(-\triangle)^{s/2}f\|_{L^{q}(\mathbb{R}^{n})}\|v(s,\cdot)\|_{L^{q'}(\mathbb{R}^{n})}ds\\
&=&\|(-\triangle)^{s/2}f\|_{L^{q}(\mathbb{R}^{n})}\int^{t}_{0}\|e^{-t(-\triangle)^{s/2}}g\|_{L^{q'}(\mathbb{R}^{n})}ds\\
&\leq&C_{n,s,q}\|(-\triangle)^{s/2}f\|_{L^{q}(\mathbb{R}^{n})}\|g\|_{L^{1}(\mathbb{R}^{n})}\int_{0}^{t}s^{-\frac{n}{sq}}ds\\
&\leq&C_{n,s,q}t^{1-\frac{n}{sq}}\|(-\triangle)^{s/2}f\|_{L^{q}(\mathbb{R}^{n})}\|g\|_{L^{1}(\mathbb{R}^{n})}
\end{eqnarray*}
for all $t>0.$ Combing the duality argument and   these two
estimates, we have
\begin{eqnarray*}
\|f\|_{L^{\infty}(\mathbb{R}^{n})}&=&\sup_{\|g\|_{L^{1}(\mathbb{R}^{n})}\leq
1, g\in \mathcal
{S}}|\langle f,g\rangle|\\
&\leq &C_{n,s,q}\left(q_{1}t^{-\frac{n}{sq_{1}}}
\left(\|f\|_{L^{p}(\mathbb{R}^{n})}+\|f\|_{Q_{\alpha}(\mathbb{R}^{n})}\right)
+t^{1-\frac{n}{sq}}\|(-\triangle)^{s/2}f\|_{L^{q}(\mathbb{R}^{n})}\right)
\end{eqnarray*}
for all $t>0$ and $p\leq q_{1}<\infty.$
 Take $$q_{1}=\log(1/t), \
 t=\left(e^{p}+\|(-\triangle)^{s/2}f\|_{L^{q}(\mathbb{R}^{n})}^{\left(1-\frac{n}{sq}\right)^{-1}}\right)^{-1}.
 $$
Then $t^{-n/(sq_{1})}=(t^{1/\log t})^{n/s}=e^{n/s}$ and
$$t^{1-\frac{n}{sq}}\|(-\triangle)^{s/2}f\|_{L^{q}(\mathbb{R}^{n})}=\left(e^{p}
+\|(-\triangle)^{s/2}f\|_{L^{q}(\mathbb{R}^{n})}^{\left(1-\frac{n}{sq}\right)^{-1}}\right)^{-(1-\frac{n}{sq})}\|(-\triangle)^{s/2}f\|_{L^{q}(\mathbb{R}^{n})}\leq
1.$$ Since we can find constant $C_{n,s,p,q}$ such that $q_{1}\leq
C_{n,s,p,q}
 \log\left(e+\|(-\triangle)^{s/2}f\|_{L^{q}(\mathbb{R}^{n})}\right),$   (\ref{mean ch Q a b 25}) holds.

\vspace{0.1in} \noindent
 {\bf{Acknowledgements.}} We would like to thank our  supervisor Professor Jie Xiao
  for suggesting the problem and kind  encouragement.

\end{document}